\input amstex.tex
\documentstyle{amsppt}
%\pagewidth{13cm}
\overfullrule=0pt
%\input amsppt.mor
%\selectlanguage{english}
\pageheight{6.9in}
\leftheadtext{C.D.Hill and M.Nacinovich}

%\rightheadtext{\titolo}

\magnification=\magstep1

%\font\tensc=cmtcsc10%\font\sevensc=cmtcsc7 \font\fivesc=cmtcsc5
\font\sc=cmcsc10
%\def\sc{\fam\scfam\tensc}%
%\textfont\scfam=\tensc

\hyphenation{pseu-do-con-ca-ve}
\hyphenation{s-o-n-d-e-r-f-o-r-s-c-h-u-n-g-s-b-e-r-e-i-c-h}
\newcount\q
\newcount\x
\newcount\t
\newcount\u

\long\def\se#1{\advance\q by 1
\x=0  \t=0 \u=0 \bigskip
\noindent
\S\number\q \quad
{\bf {#1}}\par
\nopagebreak}

%thm=teorema
%p,q parametri
\long\def\thm#1{\advance\x by 1
\bigskip\noindent%
{\sc Theorem \number\q.\number\x}
\quad{\sl #1} \smallskip\noindent}

\long\def\thml#1#2{\advance\x by 1
\bigskip\noindent
{\sc Theorem \number\q.\number\x ({#1})}
\quad{\sl #2} \smallskip\noindent}

%\lem=lemma
%parametri q,x
\long\def\lem#1{\advance\x by 1
\medskip\noindent
{\sc Lemma \number\q.\number\x}
\quad{\sl #1} \smallskip\noindent}
\long\def\prop#1{\advance\x by 1
\medskip\noindent
{\sc Proposition \number\q.\number\x}
\quad{\sl #1} \smallskip\noindent}

%\ex=esempio
%\long\def\ex{\advance\t by 1
%\medskip
%\noindent
%{\sc Example} \sn .\exn \quad}

%\os=osservazione
\long\def\os{\medskip\advance\u by 1
\noindent
{\sc Remark} \number\q.\number\u \quad}

%form=formula che si vuole etichettare
\long\def\form#1{\global\advance\t by 1
$${#1} \tag \number\q.\number\t$$}
\long\def\cor#1{\advance\x by 1
\bigskip\noindent%
{\sc Corollary \number\q.\number\x}%
\quad{\sl #1} \smallskip\noindent}

\def\dimo{\noindent{\sc Proof}\quad}

\hyphenation{pseu-do-con-vex}
\hyphenation{pseu-do-con-ca-ve}
\hyphenation{Hes-si-an}

\topmatter
\title Stein fillability and the realization of contact manifolds
\endtitle
\author C.Denson Hill and Mauro Nacinovich
\endauthor
\address C.D.Hill - Department of Mathematics, SUNY at Stony Brook,
Stony Brook NY 11794, USA \endaddress
\email dhill\@math.sunysb.edu \endemail
\address Mauro Nacinovich - Dipartimento di Matematica  -
Universit\`a di Roma ``Tor Vergata'' - 
via della Ricerca Scientifica, 1 - 00133 - ROMA, Italy \endaddress
\email nacinovi\@mat.uniroma2.it\endemail

\subjclass 35 32 53  \endsubjclass
\abstract
There is an intrinsic notion of what it means for a
contact manifold to be the smooth boundary of a Stein
manifold. The same concept has another more extrinsic 
formulation, which is often used as a convenient working
hypothesis. We give a simple proof that the two are 
equivalent. Moreover it is shown that, even though a
border always exists, it's germ is not unique; nevertheless
the germ of the Dolbeault cohomology of any border is unique.
We also point out that any Stein fillable compact contact $3$
manifold has a geometric realization in $\Bbb C^4$ via an embedding,
or in $\Bbb C^3$ via an immersion.
\endabstract
\endtopmatter
\document
Let $M$ be a smooth orientable compact real $2n+1$ dimensional manifold
without boundary
($n=1,2,3,\hdots$). Let $\Xi$ be a smooth orientable contact structure on
 $M$. The orientation of $\Xi$ is defined by a global contact form
 $\xi$ on $M$, with $\Xi=\{v\in TX\, | \, \xi(v)=0\}$, and which
is strongly non-integrable, so that
 $\omega=\xi\wedge\left(d\xi\right)^n$ is $\neq 0$ at each $x\in M$,
so defining an orientation of $M$.
We shall always take $\omega$ as the orientation of $M$, and we shall say
then that $M$ and $\Xi$ are equally oriented.\par
{\it Assume that the contact manifold $(M,\Xi)$ is the smooth boundary of a
Stein manifold
$(X,J)$.}
\par
Let us clarify this notion:
Let $X$ be a complex manifold, of dimension $(n+1)$, with a smooth
boundary $M$. 
Assuming that its complex structure $J$ is smooth up to
the boundary
$M$, it induces a smooth $CR$ structure $\left(M,HM,J_M\right)$,
$J_M:HM@>>>HM$,
 $J_M^2=-I$ of hypersurface type $(n,1)$ on $M$.
To say that a contact structure $\Xi$ on $M$ is induced by the
 $CR$ structure of $M$ means that $\Xi=HM$ are the same distribution
of $2n$-planes in $TM$. Since $M$ is a boundary, the contact structure
 $\Xi$ is orientable and a global contact form $\xi$ defines the Levi
form of $M$:
\form{ HM\ni v @>>> d\xi(J_Mv,v)\in\Bbb R\, .}
This is a Hermitian form on $HM$, for the complex structure $J_M$.
The strong non-integrability condition $\xi\wedge\left(d\xi\right)^n\neq 0$,
together with the formal integrability of the partial complex structure
 $J_M$, imply that for each $x\in M$ the Levi form
 $H_xM\ni v @>>>\Cal L(v)\in\Bbb R$ is {\sl non-degenerate}, i.e. all
its eigenvalues are different from zero. 
\par
In particular, when $M$ is the boundary of a Stein manifold $X$, the 
Levi form $\Cal L$ of $M$ is positive definite at every $x\in M$: in
this case the induced $CR$ structure is strongly pseudoconvex.
 In this situation it is customary to say that
"the contact manifold $M$ is Stein fillable by $X$".
\par
The purpose of this note is to delve into the issue of the
meaning of
the sentence in italics.\par
%%%%%%%%%%%%%%%%%%%%%%%%%%%%%%%%%%%%%%%%%%%%%%%%%%%%%%%%%%%%%%%%%%%%%%%%%%%
\se{The intrinsic notion}
Here is the issue: What is meant by saying that $M$ is the smooth boundary
of a complex
manifold $X$? If we are to enjoy the convenience and  flexibility of a
differential
topologist, and want to work in the smooth ($\Cal C^\infty$) category, then
the intrinsic
notion is clear. It goes as follows:
\roster
\item"($i$)" $\overline{X}=X\cup M$ has the structure of a $\Cal C^\infty$
manifold with
a $\Cal C^\infty$ boundary $M$, $X$ being the interior of $\overline{X}$.
\item"($ii$)" $X$ is endowed with a formally integrable almost complex
structure
$J:TX@>>>TX$, $J^2=-I$, which is $\Cal C^\infty$ up to the boundary $M$.\par
\endroster
[This much gives us a smooth induced almost-$CR$ structure
 $J_M$ on
 $M$, which in turn induces a distribution of $2n$-planes $\Xi=HM$ on $M$.
When $n=1$, there are no integrability conditions and in fact 
the $CR$ structure can be taken strictly pseudoconvex
if the corresponding contact structure is strongly non-integrable.] For Stein
fillability
we require in addition that
\roster
\item"($iii$)" $X$ is a Stein manifold.
\endroster
\os
It follows from ($ii$) via the Newlander-Nirenberg theorem that $X$ has an
atlas
of {\it interior} holomorphic coordinate charts. But it does not follow
immediately
from the above definition that $\overline{X}$ has an atlas of holomorphic
coordinate
charts [which would have to include {\it boundary charts}]. Nor does it
immediately
follow from the definition that $\overline{X}$ can be regarded as the closure
of a domain in some larger open complex manifold $\widetilde X$. See for
example the
discussion in [H1], [H2], [H3].
\se{A working hypothesis}
There has been considerable recent interest in compact contact manifolds
which are
Stein fillable, and many very 
interesting and significant results have been obtained,
especially
when $\roman{dim}_{\Bbb R}M=3$ (see e.g. [El1], [El2], [El3], [Go], [LiM]).
\par
In these articles, however, the intrinsic notion is not 
always being used; what is
being used
instead is the following convenient {\it working hypothesis}:
\roster
\item"$1^o$" The Stein manifold $X$ is an open set in a larger open
complex manifold $Y$, with $X\Subset Y$.
\item"$2^o$" There exists a real $\Cal C^\infty$ strictly plurisubharmonic
function
$\phi$ on $Y$.
\item"$3^o$" $\overline{X}=X\cup M=\{x\in Y\, | \, \phi(x)\leq 0\,\}$ with
$d\phi\left|_{\dsize M}\right. \neq 0$.
\item"$4^o$" $\phi$ is a Morse function on $Y$; i.e. $\phi$
has at most a finite number
of critical points, all of which are nondegenerate.
\endroster
This working hypothesis clearly implies the intrinsic notion, but it also
involves a
number of extrinsic elements. 
In \S 6 we give a simple proof that the intrinsic notion
is equivalent to the convenient working hypothesis.
%%%%%%%%%%%%%%%%%%%%%%%%%%%%%%%%%%%%%%%%%%%%%%%%%%%%%%%%%%%%%%%%%%%

\se{Existence and non-uniqueness of the border}
In this section we do not need that $M$ be compact, nor that $X$ be Stein.
But we will
tacitly assume that all the manifolds are paracompact (i.e. countable at
infinity).
Otherwise we place ourselves in the position of ($i$) and ($ii$) of the
intrinsic
notion.
\thm{Assume that the contact manifold $M$ is the $\Cal C^\infty$ intrinsic
boundary
of a strictly pseudoconvex complex manifold $X$. Then:
\roster
\item"(a)" $\overline{X}$ is a domain $\overline{X}\subset\widetilde X$, having
interior
$X$ and $\Cal C^\infty$ strictly pseudoconvex boundary $M$, in some open
complex
manifold $\widetilde X$.
\item"(b)" Even though a border $\widetilde X\setminus\overline{X}$ exists
by (a),
its germ
along $M$ is, in general, not unique.
\endroster}\edef\thmone{\number\q.\number\x}
\dimo (a) Since by ($i$) $\overline{X}$ is a smooth manifold with a smooth
boundary,
there is a $\Cal C^\infty$  collar, so that we can consider $\overline{X}$ as a
domain in some open real $2n+2$ dimensional smooth manifold $\Omega$. By ($ii$)
there is a complex structure tensor $J$ on $X$ which is $\Cal C^\infty$ up to
$M$, and hence induces the strictly pseudoconvex structure $J_M$ on $M$. As $J$
is assumed in ($ii$) to be $\Cal C^\infty$ up to the boundary, we may
consider its
smooth extension $\overline J$ to $\overline X$,
so $J_{M}=\overline J\left|_{\dsize HM}\right.=\overline
J\left|_{\dsize\Xi}\right. $.
Since Whitney sections over closed sets can be continued to smooth sections
over open
neighborhoods, we may, after possibly shrinking $\Omega$, extend $\overline
J$ to
a smooth {\it almost} complex structure $J_{\Omega}$ on $\Omega$, such that
${J_{\Omega}}\left|_{\dsize{\overline X}}\right.=\overline J$
satisfies the formal integrability
conditions of the Newlander-Nirenberg theorem on $\overline
X\subset\Omega$. Now
the statement (a) is the content of Theorem 1 in [HN1], where a detailed
proof is
given. It tells us that there is an open submanifold $\widetilde X$, with
$\overline X\subset\widetilde X\subset \Omega$, and a complex structure $\widetilde
J$ on
$\widetilde X$, such that $\widetilde J\left|_{\dsize\overline X}\right.=\overline J$.
The proof of that theorem involves a tricky use of Zorn's lemma, and employs an
up-to-the-boundary version of the Newlander-Nirenberg theorem, which is
valid here
since $M$ is strictly pseudoconvex (see [HJ], [Ca]).\par
This completes the proof of (a).
\smallskip
\os 
When $M$ is compact, weakly pseudoconvex
and of finite type in the sense of D'Angelo (see [DA]), the existence
of $\tilde X$ was shown by [Ch] using a much more complicated
argument. When M is compact, strictly pseudoconvex, and is
a boundary in the concrete sense (see [H1]), the existence
of $\tilde X$ was shown by [Oh] and [He]. Additional very 
interesting related results were obtained in [Le1], [Le2], [Le3].
\smallskip

(b) We give a simple counterexample to uniqueness of the germ of the border
along
$M$, even in the simple case where $X=B$ is an open ball in $\Bbb C^m$
($m=1,2,3,\hdots$) with boundary $\partial B=S^{2m-1}$. For convenience take
$B$ to be the ball of radius $\frac{1}{2}$ centered at the point
$\frac{1}{2}e_1$,
where $e_1=(1,0,\hdots,0)$. Let $D$ denote the open unit disc in $\Bbb C$,
and $\omega$
denote a suitable open neighborhood of $\overline D$, to be chosen later.
We set
$U=\omega\times D^{m-1}$ and note that $U$ is an open neighborhood of
$\overline B$
in $\Bbb C^m$. On $U$ we have the standard complex structure, which can be
described
by a single global holomorphic coordinate patch $(U;z_1,\hdots,z_m)$. We shall
construct another complex structure on $U$, also described by a single global
holomorphic coordinate patch of the form
 $(U;\widetilde\phi(z_1),z_2,\hdots,z_m)$
such that:
\roster
\item the two complex structures coincide on $\overline B$,
\endroster
while
\roster
\item"(2)" the two complex structures cannot possibly coincide on any
neighborhood in
$U$ of the point $e_1\in\partial B$.
\endroster
This means that the standard complex structure on $\overline B$ can be
extended in
inequivalent ways to the border $U\setminus\overline B$.\smallskip
Let $\alpha(z)$ denote the branch of $\sqrt{1-z}$ on $\Bbb C\setminus
[1,\infty)$
which has positive real part. On the closure $\overline D$ we define
$$\phi(z)=\cases
Az+\exp\left(-\frac{1}{\alpha(z)}\right),&\quad z\neq 1\\
A&\quad z=1\, .
\endcases$$
For every $A\in\Bbb C$ this defines a $\Cal C^\infty$ function on
$\overline D$, in
the sense of Whitney. For $|A|$ sufficiently large, it defines a
biholomorphism of
$D$ onto an open domain $G$ in $\Bbb C$. By Whitney's theorem, for large $A$,
$\phi$ extends to a smooth diffeomorphism $\widetilde\phi$ of an open neighborhood
$\omega$ of $\overline D$ in $\Bbb C$ onto a neighborhood $\Omega$ of
$\overline G$
in $\Bbb C$.\par
It follows from what was said above that the two complex structures are
equivalent
on $\overline D$, and hence on $\overline B$, yielding (1). It remains to
establish
(2): Consider the function $f(z_1,\hdots,z_m)=\widetilde\phi(z_1)$ defined on
$U$. Then
$f\left|_{B}\right.$ is holomorphic with respect to either of the two
complex structures,
and it is holomorphically extendable across $e_1$ with respect to the
second one, since
it is one of the holomorphic coordinate functions. But $f\left|_{B}\right.$
is not holomorphically extendable across $e_1$ with respect to the standard
complex
structure, because if it were extendable across $e_1$, then $Az-\phi(z)$
would have
a nonzero holomorphic extension to a neighborhood of $1$ in $\Bbb C$, while
at the
same time being flat at $1$; this gives a contradiction.
%%%%%%%%%%%%%%%%%%%%%%%%%%%%%%%%%%%%%%%%%%%%%%%%%%%%%%%%%%%%%%%%%%%%%%%%%%%
\se{Fundamental system of Stein neighborhoods}
Now we return to the situation where $M$ is compact and $X$ is Stein.
Theorem \thmone\; supplies us with an open complex manifold $\widetilde X$,
in which $\overline{X}=X\cup M$ appears as a compact domain with a smooth
strictly pseudoconvex boundary.
\thm{Assume that the compact contact manifold $M$ is the $\Cal C^\infty$
intrinsic boundary of a Stein manifold $X$. Then $\overline X$ has
a fundamental system of open Stein neighborhoods $\{Y\}$ with
 $\overline X\Subset Y\Subset\widetilde X$, for each $Y$.} 
\edef\thmtre{\number\q.\number\x}
\dimo
 This now follows from an old result
that is proved using the bumping technique of [AG],
applied to the strictly pseudoconvex domain $\overline X$ in $\widetilde X$:
by employing a finite number of small smooth bumps, one can construct an
arbitrarily small open neighborhood $Y$ of $\overline X$, such that
 $\partial Y$ is smooth and remains strictly pseudoconvex. Then using
local vanishing theorems for coherent analytic sheaves, and the Mayer-Vietoris
sequence, applied a finite number of times, it can be shown that the
restriction homomorphism
 $$r: H^q(Y,\Cal F) @>>> H^q\left(X,\Cal F|_X\right)$$
is an isomorphism for $q>0$, and any coherent analytic sheaf
 $\Cal F$
on $Y$. We have that $H^q(Y,\Cal F)\simeq H^q\left(X,\Cal F|_X\right)=0$
because $X$ is Stein. For more details, see Theorem 5 in [AH2], or
consult [AG; Propositions 16, 17, 21, 22].
%%%%%%%%%%%%%%%%%%%%%%%%%%%%%%%%%%%%%%%%%%%%%%%%%%%%%%%%%%%%%%%%%%%%%%%%%%
\se{Geometric realization of Stein fillable contact structures}
Let $n=1$, so $\roman{dim}_{\Bbb R}M=3$ and $\roman{dim}_{\Bbb C}X=2$.
\thm{Assume that the $3$-dimensional compact contact manifold $M$ is the
 $\Cal C^\infty$ intrinsic boundary of a Stein manifold $X$. Then $M$
has a smooth $CR$ embedding as a closed $CR$ submanifold of $\Bbb C^4$
(or a closed $CR$ immersion in $\Bbb C^3$).}
Note that this means that the $CR$ structure induced on $M$ from the
embedding is the same as the one $M$ inherits from being the boundary of
 $X$. In particular: {\sl the contact structure on $M$ is achieved, via
the embedding, by a complex tangent line at each point.}\par
\noindent
\dimo
Choose one of the Stein manifolds $Y\Supset \overline X$. According to the
embedding theorem for Stein manifolds (see [Bi], [Na]), $Y$ has a proper
holomorphic embedding as a closed complex submanifold of $\Bbb C^5$. The
restriction of this embedding to $M$ gives a $CR$ embedding of $M$ into
 $\Bbb C^5\subset\Bbb C\Bbb P^5$. With $N=5$ consider
 $$ M'=\{(p,r)\in M\times\Bbb C\Bbb P^N\, | \, \overline{pr}\;
\text{is tangent to $M$ at $p$}\}.$$
Then $M'$ is a smooth submanifold of $M\times\Bbb C\Bbb P^N$ of real
dimension $6$, and $M'\ni (p,r)@>>> r\in\Bbb C\Bbb P^N$ is a smooth map.
By Sard's theorem its image has measure zero in $\Bbb C\Bbb P^N$, since
 $2N>6$. By choosing a point $R_0\notin\{\text{its range}\}\cup M$,
and taking a holomorphic projection from $R_0$ to a hyperplane $\Sigma$
not containing $R_0$, we obtain a $CR$ closed immersion of $M$ into a
 $\Bbb C^{N-1}$. \par
Next consider
 $$M''=\{(p,q,r)\, | \, (p,q)\in M\times M\setminus\Delta\, ,\;
r\in\Bbb C\Bbb P^N\;\text{and}\;p,q,r\;\text{are collinear}\}.$$
Then $M''$ is a smooth manifold of real dimension $8$, and
 $M''\ni (p,q,r) @>>> r\in\Bbb C\Bbb P^N$ is a smooth map. Again by
Sard's theorem, its image has measure zero, because $2N>8$. Thus it
is possible to choose the point $R_0$ so that the $CR$ immersion
obtained above is globally one-to-one. As a result we obtain a $CR$
embedding of $M$ into $\Bbb C^4$. To obtain a $CR$ immersion into
 $\Bbb C^3$, we repeat the above projection argument with $N=4$, as
then we still have $2N>6$.
%%%%%%%%%%%%%%%%%%%%%%%%%%%%%%%%%%%%%%%%%%%%%%%%%%%%%%%%%%%%%%%%%%%%%
\os
When $n=2,3,\hdots$, so that $\roman{dim}_{\Bbb R}M\geq 5$, the result
analogous to Theorem \thmtre\; holds without any assumption of Stein 
fillability; one needs only the existence of a $CR$ structure on $M$
which is compatible with the contact structure: assume the $(2n+1)$
dimensional compact orientable contact manifold $M$ has a smooth
 $CR$ structure of type $(n,1)$ which induces the given contact structure
and is strictly pseudoconvex. By a theorem of 
Boutet de Monvel [BM], $M$ has a smooth $CR$ embedding into 
 $\Bbb C\Bbb P^N$, for some $N$. Then we can repeat the
argument above, and obtain that $M$ has a $CR$ embedding into
 $\Bbb C^{2n+2}$,
or a $CR$ immersion into $\Bbb C^{2n+1}$.
The contact structure on $M$ is then achieved, via the embedding, by a
tangent affine $\Bbb C^n$ at each point.
\par
For $CR$ manifolds which are not of hypersurface type,
see [HN2].
%%%%%%%%%%%%%%%%%%%%%%%%%%%%%%%%%%%%%%%%%%%%%%%%%%%%%%%%%%%%%%%%%%%%%%%%%
\se{Equivalence of the intrinsic notion and the working hypothesis}
We return to the situation of \S 1 and \S 2. 
\thm{Assume that the compact contact manifold $M$ is the $\Cal C^\infty$
intrinsic boundary of a Stein manifold $X$. Then the working hypothesis
$1^o$, $2^o$, $3^o$, $4^o$ of \S 2 are satisfied.}
\dimo Since a Stein manifold $Y$ has an exhaustion by a smooth strictly
plurisubharmonic function, we obtain $1^o$ and $2^o$ from 
Theorems \thmone\; and \thmtre. To demonstrate $3^o$ we proceed as follows:
Fix a Stein neighborhood $Y$ of $\overline X$ in $\widetilde X$, a strictly
plurisubharmonic function $\psi$ on $Y$, and a Hermitian metric on $Y$.
As $X$ is a domain in $Y$, there exists a global defining function
 $\rho\in\Cal C^\infty(Y)$ such that:
 $$\overline X=\{x\in Y\, | \, \rho(x)\leq 0\}\, ,\qquad
 d\rho|_M\neq 0\, ,$$
(see [AH1; Proposition 1.1]). Since $M$ is strictly pseudoconvex the 
Levi form $\Cal L(\rho)$ is positive definite at each point of $M$; 
i.e. has $n$ positive eigenvalues. To obtain $(n+1)$ positive eigenvalues
for the complex Hessian $i\partial\bar\partial\rho$ near $M$, we replace
 $\rho$ by a modified global defining function
 $$\widetilde\rho=\frac{1}{\lambda}\{e^{\lambda\rho}-1\}\, ,$$
with the constant $\lambda>0$ chosen sufficiently large. It is easy to
verify that there is an open neighborhood $U$ of $M$ in $Y$ 
in which $\widetilde\rho$ is strictly plurisubharmonic, and
 $d\widetilde\rho\neq 0$. Next we modify $\widetilde\rho$ to make it strictly
plurisubharmonic in an open neighborhood $V$ of $\overline X$, and
establish $3^o$: Let $\chi(\rho)$ be a smooth real convex function of
the real variable $\rho$, such that $\chi(\rho)=\rho$ for
 $\rho\geq -\delta$ and $\chi(\rho)=-2\delta$ for $\rho\leq -3\delta$,
where $\delta>0$ is chosen so small that 
 $\{x\in\overline X\, | \, -3\delta\leq\widetilde\rho(x)\leq 0\}\subset U$.
Let  $K\subset X$ be a compact set such that
 $K\supset\{x\in X\, | \,\widetilde\rho(x)\leq-\delta\}$. Choose a nonnegative
smooth cutoff function $\mu\in\Cal C^\infty_0(X)$ such that
 $\mu=1$ on a neighborhood of $K$. Consider the function:
 $$\phi=\chi(\widetilde\rho)+\epsilon\mu\psi\, ,$$
with a small constant $\epsilon>0$. Then
 $d\phi|_M=d\widetilde\phi|_M=d\rho|_M\neq 0$ and
 $\overline X=\{x\in Y \, | \, \phi(x)\leq 0\,\}$ for $\epsilon>0$
taken sufficiently small. The function $\chi(\widetilde\rho)$ is smooth and
weakly plurisubharmonic on $V=X\cup U$. The function $\phi$ is strictly
plurisubharmonic in $V$ for sufficiently small $\epsilon>0$. This
establishes $3^o$ without destroying $2^o$.\par
The function $\psi$ can be chosen at the beginning to be a Morse function
on $Y$; see [AF]. Hence by construction there is an $\eta>0$ such that
 $\phi$ has no critical points on
 $\{x\in V\, | \, -\eta\leq\widetilde\rho(x)\leq\eta\,\}$,
and at most only a finite number of nondegenerate critical points for
 $\{x\in V\, | \, \widetilde\rho(x)\leq -3\delta\,\}$. To obtain
 $4^o$, we need to eliminate any degenerate critical points of $\phi$ in
 $\{x\in V\, |\, -3\delta<\widetilde\rho(x)<-\eta\}$. 
Let $\nu\in\Cal C^\infty_0(X)$, $0\leq\nu(x)\leq 1$, be a smooth
cutoff function with $\nu=1$ on the set $\{x\, | \, \widetilde\rho(x)\leq-\eta\}$.
By Sard's theorem we can approximate $\phi$, in the $\Cal C^2$-norm
on any compact subset of $V$, by a smooth function $\widetilde\phi$ which has
only nondegenerate critical points; hence $\widetilde\phi$ remains strictly
plurisubharmonic. Set
  $$\phi_1=\nu\widetilde\phi+(1-\nu)\phi\, .$$
Then for $\phi_1-\widetilde\phi=(1-\nu)(\phi-\widetilde\phi)$
there is an estimate
  $$\left|\phi_1-\widetilde\phi\right|_2\leq \text{const}\,
\left|\phi-\widetilde\phi\right|_2\, ,$$
where the norms are $\Cal C^2$-norms taken over some compact subset
 $L\Subset V$, with $\overline X\subset\overset{o}\to{L}$. So by taking
a sufficiently good approximation $\widetilde\phi$ to $\phi$, the function
 $\phi_1$ satisfies $1^o$, $2^o$, $3^o$, $4^o$; hence the proof is complete.
%%%%%%%%%%%%%%%%%%%%%%%%%%%%%%%%%%%%%%%%%%%%%%%%%%%%%%%%%%%%%%%%%%%%%%%
\se{Cohomology of the border}
In spite of the fact that the germ of the border
 $\widetilde X\setminus\overline X$
is not unique, it turns out that the germ of its Dolbeault cohomology
is unique:
\thm{Assume that the compact contact manifold $M$ is the $\Cal C^\infty$
intrinsic boundary of a Stein manifold $X$. Then for any choice of the 
 $\widetilde X$, in which $\overline X$ is a domain, and for any choice of
the Stein neighborhood $Y$, $\overline X\subset Y\Subset\widetilde X$, and
for any $0\leq p\leq n+1$, we have:
\roster
\item $H^{p,q}(Y\setminus X)\simeq H^{p,q}(M)=0$ for $0<q<n$,
\item $H^{p,n}(Y\setminus X)\simeq H^{p,n}(M)$,
\item $H^{p,n+1}(Y\setminus X)=0$,
\endroster
with 
\roster
\item"(4)"
 $\roman{dim}_{\Bbb C}H^{p,n}(Y\setminus X)=\infty.$\endroster}
Here $H^{p,q}(Y\setminus X)$ denotes the Dolbeault cohomology of smooth
 $\bar\partial$-closed $(p,q)$-forms on $Y\setminus X$ modulo those
which are $\bar\partial$ exact in $Y\setminus X$. Note that
 $Y\setminus X=\left(Y\setminus\overline X\right)\cup M$ has smooth
boundary $M$, and we are requiring here that the differential forms
be $\Cal C^\infty$ up to $M$. $H^{p,q}(M)$ denotes the
 $\bar\partial_M$-cohomology of tangential $\bar\partial_M$-closed
smooth $(p,q)$-forms on $M$, modulo those that are $\bar\partial_M$-exact
on $M$. 
\par
The results (1), (2), (3), (4) are direct consequences of [AH1], [AH2]; see
Theorems 5 and 7], or see Theorem 7.2 in [HN3], and [La].
\os 
When $q=0$ and $0\leq p\leq n+1$ we have that
 $H^{p,0}(Y\setminus X)\simeq H^{p,0}(Y)$ and
 $H^{p,0}\left(\overline X\right)\simeq H^{p,0}(M)$, 
see [AH1].
\par

\enddocument